\documentclass[12pt,twoside]{amsart}

\newtheorem{Thm}{Theorem}
\newtheorem{Cor}[Thm]{Corollary}
\newtheorem{Lemma}[Thm]{Lemma}
\newtheorem{Prop}[Thm]{Proposition}
\newtheorem*{Prop*}{Proposition}
\newtheorem{Conj}[Thm]{Conjecture}

\theoremstyle{definition}

\newtheorem{Ex}{Example}

\input{diagrams}

\title{Fredholm modules over certain group C*-algebras}
\author[T.~Hadfield]{Tom Hadfield}
\address{Department of Mathematics, University of California, Berkeley, CA 94720}
\email{hadfield@math.berkeley.edu}
\date{\today}

\begin{document}

\begin{abstract}
Motivated by the search for new examples of ``noncommutative manifolds'', we study the noncommutative geometry of the group C*-algebras of various discrete groups.
 The examples we consider are the infinite dihedral group ${\bf Z} {\times_{\sigma}} {\bf Z}_2$ and the semidirect product  group ${\bf Z} {\times_{\sigma}} {\bf Z}$.
 We present a unified treatment of the K-homology and cyclic cohomology of these algebras.
\end{abstract}

\maketitle

\section{Introduction}

Recently, there has been a great deal of interest in the notion of a ``noncommutative manifold''. 
This problem has been studied from several different viewpoints, and a consensus on a final definition seems to be quite far off. 
 However, the study of examples that one would surely want such a  definition to encompass has been generating very interesting new mathematics.  

In \cite{connes96} Connes formulated the notion of a ``noncommutative Riemannian spin manifold'', within his functional-analytic  framework of noncommutative geometry.
 For a given C*-algebra, 
this consists of a spectral triple \cite{connes94} over the algebra 
 (a *-representation of the C*-algebra as bounded operators on a Hilbert space, together with an unbounded ``Dirac'' operator)
 with additional structures (such as a reality operator) satisfying an appropriate list of axioms. 
In the situation where the algebra is commutative, it is a nontrivial result \cite{fgv}, \cite{rennie} that the algebra must necessarily consist of differentiable functions on a Riemannian spin manifold, with the manifold being recovered as the spectrum of the algebra. 
 Until very recently the list of examples of noncommutative algebras which could be equipped with the structure of a noncommutative geometry was rather short. 
 The most useful example was given by the noncommutative tori \cite{rieffel81}.

The approach via quantum groups has been much more examples-driven : see for example  Majid \cite{majid00}, Schmudgeon \cite{schmudgeon}, Woronowicz \cite{woro}.
 Rather than taking as a starting point an axiomatic framework and then looking around for examples that fit into it, 
 they began with a case by case study of 
 the many interesting and provocative examples coming from quantum groups, 
and for each one 
 try to see which of the structures of ordinary commutative geometry can be translated meaningfully into this picture. 
A particular focus of this approach is to classify the possible differential calculi over a given algebra.

An instructive example illustrating the differences between the two approaches is the algebra $C( {SU_q}(2))$, ``continuous functions on the quantum $SU(2)$''. 
 This algebra was thoroughly investigated using the tools of Connes' noncommutative geometry by Masuda, Nakagame and Watanabe in \cite{mnw90}.
 In particular the K-theory and K-homology were calculated, and the generators exhibited : all four groups are in fact isomorphic to ${\bf Z}$. 
 So from this viewpoint the algebra appeared to be quite simple. 

Woronowicz's original work \cite{woro} classified the covariant differential calculi over\\ 
 $C( {SU_q}(2))$. 
 Subsequently, Schmudgeon \cite{schmudgeon}
 showed that there are no representations of these calculi by bounded commutators, a result which seems to be in direct opposition to Connes' approach.

In the last year there has been a flurry of interest in producing more examples of noncommutative spectral triples (see \cite{cl} and many others). 
 This work culminated in Varilly's paper \cite{varilly01}, in which Rieffel's deformation quantization techniques \cite{rieffel93} were used to show that the noncommutative spheres of Connes and Landi are quantum homogeneous spaces for certain compact quantum groups,  
and also the work of Connes and Dubois-Violette \cite{cdv}, which gives a complete description from K-theoretic considerations of three-dimensional noncommutative spherical manifolds. 

The approach we take in this paper is as follows. 
 We are motivated by the search for further examples of ``noncommutative manifolds'', in the sense of Connes. 
 Rather than beginning with well-known geometries over commutative C*-algebras, and then deforming the product to give noncommutative algebras as in \cite{cl},
 we instead take as our starting point the group C*-algebras of various discrete groups and investigate the extent to which they can be furnished with the structures of noncommutative geometries. 
 In general there need not be a unique way to do this. 
 
The examples we shall consider in this paper are the group C*-algebras of the infinite dihedral group ${\bf Z} {\times_{\sigma}} {\bf Z}_2$ and  the semidirect product group ${\bf Z} {\times_{\sigma}} {\bf Z}$.
 Since a noncommutative geometry over a C*-algebra $A$ consists of a spectral triple over $A$ (together with additional structure), and the bounded formulation of spectral triples are Fredholm modules, which make up the Kasparov K-homology groups ${KK^i}(A,{\bf C})$, $(i=0,1)$, 
 we start by  calculating  the K-homology groups of the algebras under investigation and look for the generating Fredholm modules.
 Having exhibited the Fredholm modules, we demonstrate their nontriviality and linear independence from one another by calculating the pairings of their Chern characters with the generators of K-theory. 
 For the infinite dihedral group, we apply Burghelea's results \cite{bu85} to calculate the cyclic cohomology of the group ring, and then via a seperate direct calculation exhibit the 0-cocycles that generate the periodic even cyclic cohomology.  

 We note that although it is well-known \cite{connes94} that corresponding to every 
 (theta summable) Fredholm module there is a corresponding spectral triple, it is not clear that every spectral triple can be equipped with the additional structure of a noncommutative geometry. 
 Furthermore, several of the Fredholm modules we exhibit are 
 built from non-faithful representations 
and would thus correspond to noncommutative geometries of dimension smaller than that of the algebra. 
 The problem of constructing the corresponding spectral triples and furnishing them with the additional structures of noncommutative geometries will be the subject of a later paper.

\section{Fredholm modules as K-homology}

We begin with some general preliminaries about Fredholm modules and K-homology.  Recall that a Fredholm module over a *-algebra 
$A$ 
is a triple 
$(H,\pi,F)$,
where $\pi$ is a *-representation of 
$A$ 
as bounded operators on the Hilbert space 
$H$. 
The operator 
$F$ is a selfadjoint element of  
${\bf B}(H)$,
satisfying $F^2 =1$, 
such that the commutators 
$[F,\pi(a)]$ are compact operators for all 
$a \in A$. 
Such a Fredholm module is called odd.

An even Fredholm module is the above data, 
together with a 
${\bf Z}_2$-grading 
of the Hilbert space $H$,
given by a grading operator 
$\gamma \in {\bf B}(H)$
with 
$\gamma = {\gamma}^{*}$, 
$\gamma^2 =1$,
$[\gamma, \pi(a)]=0$ 
for all 
$a \in A$, 
and
$F\gamma = - \gamma F$. 
In general the *-algebra 
$A$ 
will be a dense subalgebra of a C*-algebra, 
closed under holomorphic functional calculus.
Fredholm modules should be thought of as abstract elliptic operators, 
since they are motivated by axiomatizing the important properties of elliptic pseudodifferential operators on closed manifolds. 

This definition is due to Connes \cite{connes94}, p288.
In Kasparov's framework the K-homology groups are given by specialising the second variable in the KK-functor to be 
the complex numbers ${\bf C}$.
Equivalence classes of
even Fredholm modules make up the even K-homology group ${KK^0}(A,{\bf C})$.
Odd Fredholm modules make up the odd K-homology ${KK^1}(A,{\bf C})$. 
 A Fredholm module is said to be degenerate if 
 $[F, \pi(a)] =0$ for all $a \in A$. 
 Degenerate Fredholm modules represent the identity element of the corresponding K-homology group.

Two simple examples of an even and an odd Fredholm module, that we will use extensively in the sequel, are as follows:

\begin{Ex}
\label{evenfred}
Given a C*-algebra $A$, with a *-homomorphism $\phi : A \rightarrow {\bf C}$, 
we construct a canonical even Fredholm module
 ${\bf z}_0 \in {KK^0}(A, {\bf C})$ : 
\begin{equation}
{\bf z}_0 = (
H_0 = {\bf C}^2, \pi_0 = \phi \oplus 0 , 
F_0 = \left(
\begin{array}{cc}
0 & 1 \cr
1 & 0 \cr
\end{array}
\right),
\gamma = 
\left(
\begin{array}{cc}
1 & 0 \cr
0 & -1 \cr
\end{array}
\right)
).
\end{equation}
 In general ${\bf z}_0$ may well represent a trivial element of the even K-homology of $A$ (for example, if $\phi$ is the zero homomorphism.) 
However, if $A$ is unital, and $\phi$ is a nonzero  *-homomorphism, the Chern character of ${\bf z}_0$ pairs nontrivially with $[1] \in {K_0}(A)$, 
 showing that ${\bf z}_0$ is a  nontrivial element of K-homology, 
and also that $[1] \neq 0 \in {K_0}(A)$. More precisely :

\begin{Lemma} If $A$ is unital, and $\phi$ nonzero, then 
 $< {ch}_{*}( {\bf z}_0 ), [1]> =1$.
\end{Lemma}
\begin{proof} 
 Here, 
${ch}_{*} : {KK^0}(A,{\bf C}) \rightarrow {HC^{even}}(A)$, 
is the even Chern character as defined in \cite{connes94}, p295, mapping the even K-homology of $A$ into even periodic cyclic cohomology,
and $<.,.>$ denotes the pairing between K-theory and periodic cyclic cohomology defined in \cite{connes94}, p224.
 We have
\begin{equation}
 < {ch_{*}}( {\bf z}_0 ), [1] > = 
{{\lim}_{n \rightarrow \infty} } {(n!)}^{-1} \psi_{2n}  (1,...,1),
\end{equation}
 where (for each $n$) $\psi_{2n}$ is the cyclic $2n$-cocycle defined by
\begin{equation}
 \psi_{2n} ( a_0, a_1, ..., a_{2n}) = 
(-1)^{n(2n-1)} \Gamma (n+1) 
 Tr( \gamma \pi_0 (a_0) [ F_0 , \pi_0 (a_1)] ... [F_0 , \pi_0 (a_{2n})]). 
\end{equation}
 Since $\Gamma(n+1) = n!$ it follows that
\begin{equation}
< {ch_{*}}( {\bf z}_0 ), [1] > = 
{{\lim}_{n \rightarrow \infty} }
 {(-1)}^n
 Tr( \gamma \pi_0 (1) {[F_0 , \pi_0 (1)]}^{2n}).
\end{equation} 
 Now,  
$[F_0 , \pi_0 (1) ] = 
\left(
\begin{array}{cc}
0 & -1 \cr
1 & 0 \cr
\end{array}
\right),$
 hence
 $\gamma \pi_0 (1) {[F_0 , \pi_0 (1)]}^{2n} =$
$(-1)^n 
\left(
\begin{array}{cc}
1 & 0 \cr
0 & 0 \cr
\end{array}
\right)$.
 Therefore
\begin{equation} 
< {ch_{*}}( {\bf z}_0 ), [1] > = 
 {{\lim}_{n \rightarrow \infty} } {(-1)}^n 
Tr( 
(-1)^n 
\left(
\begin{array}{cc}
1 & 0 \cr
0 & 0 \cr
\end{array}
\right)) = 1
\end{equation}
as claimed. 
\end{proof}
\end{Ex}

\begin{Ex}
\label{oddfred}
 Let $A$ be a C*-algebra, together with a *-homomorphism 
$\phi : A \rightarrow {\bf C}$ and an automorphism $\alpha$ implementing an action of ${\bf Z}$. 
We describe a canonical odd Fredholm module
 ${\bf z}_1 \in {KK^1}( A \times_{\alpha} {\bf Z}, {\bf C})$ : 
\begin{equation}
 {\bf z}_1 = (H_1 = {l^2}( {\bf Z}), \pi_1 , F_1)
\end{equation}
 Take 
$\pi_1 : A \times_{\alpha} {\bf Z} \rightarrow {\bf B}( {l^2}({\bf Z}))$ 
to be defined by 
 \begin{equation}
(\pi_1 (a) \xi) (n) = \phi( {\alpha^{-n}}(a)) \xi(n), \quad  
(\pi_1 (V) \xi) (n) = \xi(n-1),
 \end{equation} 
 for $\xi \in {l^2}({\bf Z})$, $a \in A$, and $V$ the unitary implementing the action of ${\bf Z}$ on $A$ 
(via $Va{V^{*}} = \alpha (a)$).  
 Then $\pi_1$ is the usual representation of $A \times_{\alpha} {\bf Z}$ induced from the representation $\phi$ of $A$.
 We take 
\begin{equation}
F_1 \xi(n) = sign(n) \xi(n) =
\left\{
\begin{array}{cc}
 \xi(n) & : n \geq 0 \cr
- \xi(n) & : n<0 \cr
\end{array}
\right. 
\end{equation}
 It is immediate that $[F_1 , \pi_1 (a)] =0$ for all $a \in A$, and that  
  $[F_1 , \pi_1(V)]$ is a rank-one operator and hence compact. 
Nontriviality of ${\bf z}_1$ (even if $\phi$ is the zero homomorphism) follows from :

\begin{Lemma}
$< {ch}_{*} ({\bf z}_1 ), [V]> =1$. 
\end{Lemma}
\begin{proof}
 Again, 
${ch}_{*} : {KK^1}(A,{\bf C}) \rightarrow {HC^{odd}}(A)$, 
is the odd Chern character as defined in \cite{connes94}, p296, mapping the odd K-homology of $A$ into odd periodic cyclic cohomology,
and $<.,.>$ denotes the pairing between K-theory and periodic cyclic cohomology \cite{connes94}, p224. It is straightforward to calculate this pairing directly, as in the previous example, but it is quicker to appeal to Connes' index theorem \cite{connes94}, p296, which states that 
\begin{equation}
< {ch_{*}}( {\bf z}_1 ), [V]> 
= Index(EVE)
\end{equation}
 where $E = {\frac{1}{2}}(1+F)$ is the natural orthogonal projection 
 ${l^2}({\bf z}) \rightarrow {l^2}({\bf N})$. We have 
\begin{equation}
Index(EVE) = dim ker(EVE) - dim ker(E{V^{*}}E) = 1 -0 =1,
\end{equation}
 hence the result.
  This shows that ${\bf z}_1$ is a nontrivial Fredholm module, and also that $[V] \neq 0 \in$ 
${K_1}(A \times_{\alpha} {\bf Z})$. 
\end{proof}

 These Fredholm modules ${\bf z}_0$ and ${\bf z}_1$ can both be defined more generally (for example by taking instead 
$\phi : A \rightarrow {\bf B}(H)$, where $H$ is a finite dimensional Hilbert space) 
but the above formulation will be sufficient for our purposes.  
\end{Ex}

A useful tool for calculating the K-homology groups of C*-algebras is the following corollary of the universal coefficient theorem of Rosenberg and Schochet
(see Blackadar \cite{blackadar}, p234).

\begin{Prop}
\cite{rs}
\label{kgpsfreeab}
Let $A$ be a separable C*-algebra,
belonging to the  ``bootstrap class".
If the K-groups 
${K_{i}}(A)$ ($i=0,1$) 
are free abelian, then so are the K-homology groups. 
In fact  
${KK^{i}}(A, {\bf C}) \cong {K_{i}}(A)$ (as abelian groups).
\end{Prop}

We now consider the six term cyclic exact sequence for K-homology of crossed products by ${\bf Z}$, dual to the Pimsner-Voiculescu sequence for K-theory, as described in \cite{blackadar}, p199.

Recall \cite{pv} that associated to any crossed product algebra $A \times_{\alpha} {\bf Z}$  
is the following semisplit short exact sequence of C*-algebras, 
the Pimsner-Voiculescu ``Toeplitz extension" 
\begin{equation}
\label{toeplitzext}
 0 \rightarrow {A \otimes {\bf K}} \rightarrow {T_{\alpha}} \rightarrow  {A \times_{\alpha} {\bf Z}} \rightarrow  0.
\end{equation}
 Here ${T_{\alpha}}$ 
is the C*-subalgebra of 
$( {A {\times}_{\alpha} {\bf Z}}) \otimes {\it T} $
generated by 
$a \otimes 1$, 
$a \in A$ 
and 
$ V \otimes f$,
where 
$V$ 
is the unitary implementing the action of 
$\alpha$ 
on 
$A$, 
and 
$f$ 
is the non-unitary isometry generating the ordinary Toeplitz algebra
 $T$, 
that is  
$f \in {\bf B}({l^2}({\bf N}))$,
$f{e_n} = {e_{n+1}}$.
 This extension defines the Toeplitz element 
${\bf x} \in {KK^1}( {A \times_{\alpha} {\bf Z}}, A) $.  

Applying the K-functor gives the Pimsner-Voiculescu six term cyclic sequence for K-theory. 
The corresponding six term cyclic sequence for K-homology is:

\begin{diagram}
{{KK^0}( A  , {\bf C})} & \lTo^{id - {\alpha}^{*}} & {{KK^0}( A, {\bf C})} & \lTo^{i^{*}} & {{KK^0}(A {\times_{\alpha}} {\bf Z}, {\bf C})}  \\
\dTo^{\partial_0} & & & & \uTo^{\partial_1} \\
{{KK^1}( A {\times_{\alpha}} {\bf Z} , {\bf C})} & \rTo^{i^{*}} & {{KK^1}( A, {\bf C})} & \rTo^{id - {\alpha}^{*}} & {{KK^1}( A, {\bf C} )}\\
\end{diagram}

Here $i$ denotes the canonical inclusion map $i : A \hookrightarrow A \times_{\alpha} {\bf Z}$. 
The vertical maps 
${\partial_0}$ and  ${\partial_1}$  
are given by taking the Kasparov product with the  Toeplitz element:
\begin{equation}
\partial_i : {KK^i}(A,{\bf C}) \rightarrow {KK^{i+1}}( A \times_{\alpha} {\bf Z}, {\bf C}), \quad
{\bf z} \mapsto {\bf x} {\hat{\otimes}}_A {\bf z}
\end{equation}
 These
 morphisms are studied in detail in \cite{hadfield}.   
 This sequence formulated in terms of Ext appears in the original paper of Pimsner and Voiculescu \cite{pv}.
However, the relationship between Ext and the Fredholm module picture of K-homology is not transparent.

Let $A$ be unital C*-algebra, and $\phi : A \rightarrow {\bf C}$ a nonzero *-homomorphism. 
Then the  Fredholm modules ${\bf z}_0$ and ${\bf z}_1$ described above (Examples \ref{evenfred} and \ref{oddfred}) are related via the morphism  $\partial_0$ as follows.

\begin{Prop}
\cite{hadfield}
\label{partialzero}
Under the map ${\partial_0}$ 
we have $\partial_0 ( {\bf z}_0 ) = {\bf z}_1$. 
\end{Prop}

\section{The infinite dihedral group}

We give a unified treatment of the noncommutative geometry of  the group C*-algebras of the infinite dihedral group 
$\Gamma =$ ${\bf Z} \times_{\sigma} {\bf Z}_2$ 
and the semidirect product group $G=$ ${\bf Z} \times_{\sigma} {\bf Z}$, where in each case the action $\sigma$ on ${\bf Z}$ is by inversion. 
 The infinite dihedral group
$\Gamma$ is also a free product 
$ {{\bf Z}_2} * {{\bf Z}_2}$. 
We consider $G$ to be the underlying set ${\bf Z}^2$, with multiplication 
$(m,n)*(p,q) = (m+ {(-1)^n} p , n + q)$. 
These two groups are related as follows. 
We have a short exact sequence of groups 
\begin{equation}
0 \rightarrow {\bf Z} \rightarrow^{i} {\bf Z} \times_{\sigma} {\bf Z} 
\rightarrow^{q} {\bf Z} \times_{\sigma} {\bf Z}_2 \rightarrow 0.
\end{equation}
 The maps $i$ and $q$ are given by 
$i(n) = (0,2n)$, and $q(m,n) = (m, (-1)^n)$. 
 The left hand group is the centre $Z(G)$ of $G$, 
\begin{equation}
Z(G) = \{ (0,2n) : n \in {\bf Z} \} \cong {\bf Z},
\end{equation}
 whereas the righthand term is the infinite dihedral group.
 Note that $G$ is torsion-free. 

We first consider the group C*-algebra of the infinite dihedral group $\Gamma$. 
The K-theory of this algebra is well known. 
We calculate the K-homology, exhibit the generating Fredholm modules and calculate their pairings with the generators of K-theory. 
 We calculate the automorphism group $Aut(\Gamma)$ of the group $\Gamma$, and its action on the K-theory and K-homology of the group C*-algebra. 
 Finally, we calculate the cyclic cohomology of the group ring ${\bf C}\Gamma$, and of an appropriate dense ``smooth'' subalgebra.

\section{K-theory and K-homology}

Since $\Gamma$ is a semidirect product (and hence amenable) 
${\bf Z} {\times_{\sigma}} {\bf Z}_2$
 the group C*-algebra 
$A={C^{*}}(\Gamma) \cong {C^{*}_{r}}(\Gamma)$ 
 is isomorphic to 
${C({\bf T})} {\times_{\sigma}} {\bf Z}_2$,
where the 
${\bf Z}_2$-action 
$\sigma$ 
is given by
$( \sigma f)(z) = f({\bar z})$, 
for 
$z \in {\bf T}$,
 $f \in C({\bf T})$.
The algebra $A$ is generated by unitaries 
$e$ and 
$S$ 
satisfying
\begin{equation}
e = {e^{*}}, \quad
{e^2} = 1, \quad
eSe= {S^{*}} = {S^{-1}}.
\end{equation}
The K-theory of $A$ is well known. 
\begin{Prop}
\label{kthyinfdih}
\cite{blackadar}, \cite{kumjian}, \cite{lance}. 
We have 
${K_0}(A) \cong {\bf Z}^3$, 
${K_1}(A) = 0$
with the generators of 
${K_0}(A)$ 
being given by the equivalence classes of the projections  
$1$, 
${\frac{1}{2}}(1 + e)$ and  
${\frac{1}{2}}(1 + eS)$ 
in 
$A$.
\end{Prop} 

\begin{Thm}
\label{khominfdih}
The K-homology of $A$ is given by 
${KK^0}(A,{\bf C}) \cong {\bf Z}^3$, 
${KK^1}(A,{\bf C}) \cong 0$.
A basis for the even K-homology is given by Fredholm modules 
${\bf w}_0$, ${\bf w}_1$ and  ${\bf w}_2$ which we describe below. 
\end{Thm} 
 \begin{proof} 
We use 
the universal coefficient theorem 
to calculate the K-homology.
Since the K-groups are free abelian, 
it follows from Proposition \ref{kgpsfreeab} that 
 ${KK^{0}}(A, {\bf C}) \cong {\bf Z}^3$, 
${KK^{1}}(A, {\bf C}) \cong 0$, 
as claimed. 
 We will exhibit the generating even Fredholm modules. 
 Since $A \cong C({\bf T}) \times_{\sigma} {\bf Z}_2$,  
 we start by noting the K-homology of $C({\bf T})$.

\begin{Lemma}
\label{onetorus}
Both the even and odd K-homology of $C({\bf T})$ are isomorphic to ${\bf Z}$. 
 We exhibit the generating Fredholm modules. 
\end{Lemma}
\begin{proof} This follows from Prop \ref{kgpsfreeab}, since the K-groups of $C({\bf T})$ are both isomorphic to ${\bf Z}$. 

The generator of ${KK^0}(C({\bf T}),{\bf C})$ is the even canonical Fredholm module ${\bf z}_0$ (Example \ref{evenfred}) corresponding to the *-homomorphism $\phi : C({\bf T}) \rightarrow {\bf C}$, 
$U \mapsto 1$. Explicitly,
\begin{equation}
{\bf z}_0 = ( {\bf C}^2, \phi \oplus 0, 
F= \left(
\begin{array}{cc}
0 & 1 \cr
1 & 0 \cr
\end{array}
\right),
\gamma =
\left(
\begin{array}{cc}
1 & 0 \cr
0 & -1 \cr
\end{array}
\right)
).
\end{equation}
 Since $C({\bf T})$ is a crossed product ${\bf C} \times_{id} {\bf Z}$ (via the trivial action)  the generator of 
${KK^1}(C({\bf T}),{\bf C})$ 
is the odd Fredholm module 
${\bf z}_1 = (H_1 ={l^2}({\bf Z}), \pi_1, F_1)$, (Example \ref{oddfred}) 
where
$\pi_1(U) e_n = e_{n+1}$, 
and 
$F_1 {e_n} = sign (n) e_n$.  

We note that under the Baum-Connes assembly map \cite{bch}
\begin{equation}
 \mu : {KK^i}( C_0 ( B{\bf Z}), {\bf C}) \cong {KK^i}( C({\bf T}),{\bf C})
 \rightarrow 
{K_i}(C({\bf T}))
\end{equation}
 we have ${\bf z}_0 \mapsto \pm [1]$, 
${\bf z}_1 \mapsto \pm[U]$.
\end{proof}

We modify these Fredholm modules to give an even Fredholm module ${\bf w}_0$ over $A$, modelled on the corresponding ${\bf z}_0$ for $C({\bf T})$, and two even Fredholm modules ${\bf w}_1$ and ${\bf w}_2$ over $A$ induced from the odd Fredholm module ${\bf z}_1$ over $C({\bf T})$. 

The first generator of the even K-homology is given by the canonical even Fredholm module ${\bf w}_0$ (Example \ref{evenfred}) 
corresponding to the *-homomorphism   
$\phi : A \rightarrow {\bf C}$,
defined on generators by
$\phi : S,e \mapsto 1$.
 Then 
\begin{equation}
{\bf w}_0 = ( H_0 ={\bf C}^2 , \pi_0 = \phi \oplus 0, F_0 = \left(
\begin{array}{cc}
0 & 1 \cr
1 & 0 \cr
\end{array}
\right), 
\gamma=
\left(
\begin{array}{cc}
1 & 0 \cr
0 & -1 \cr
\end{array}
\right)  
).
\end{equation}
Let 
$P$ 
be any of the projections 
$1$, 
${\frac{1}{2}}(1+e)$, 
${\frac{1}{2}}(1+eS)$.
Then it is immediate that
$\pi_0 (P) =
\left(
\begin{array}{cc}
1 & 0 \cr
0 & 0 \cr
\end{array}
\right)$, 
and 
\begin{equation}
\gamma \pi_0 (P) {[F_0 , \pi_0 (P)]}^{2k}=
{(-1)}^{k}
\left(
\begin{array}{cc}
1 & 0 \cr
0 & 0 \cr
\end{array}
\right).
\end{equation}
 Hence 
$<{ch_{*}}({\bf w}_0), [P]> =1$, 
 for each of $P=$ $1$, $P_1$, $P_2$.

Now we define the Fredholm module ${\bf w}_1 \in {KK^0}(A,{\bf C})$.
First of all, consider 
$A$ 
acting on 
${l^2}({\bf Z})$,
 with orthonormal basis 
${\{ {e_{n}} \}}_{n \in {\bf Z}}$, 
via 
\begin{equation}
S{e_n} ={e_{n+1}}, \quad 
e{e_{n}} = {e_{-n}}
\end{equation}
Define the diagonal operator 
$F$ 
by
$F {e_{n}}=sign(n) {e_{n}}$. 
Let 
$H = {l^2}({\bf Z}) \oplus {l^2}({\bf Z})$.
We define a *-homomorphism 
${\pi_1} : A \rightarrow {\bf B}(H)$ 
by

\begin{equation}
\label{pione}
{\pi_1}(S)=
\left(
\begin{array}{cc}
S & 0 \cr
0 & S \cr
\end{array}
\right), \quad
{\pi_1}(e) =
\left(
\begin{array}{cc}
e & 0 \cr
0 & -e \cr
\end{array}
\right).
\end{equation}
Define
\begin{equation}
F_1 =
\left(
\begin{array}{cc}
0 & i{F} \cr
-i{F} & 0 \cr
\end{array}
\right), 
\quad
\gamma = 
\left(
\begin{array}{cc}
1 & 0 \cr
0 & -1 \cr
\end{array}
\right).
\end{equation}

It is easy to verify that $[F_1 ,{\pi_1}(S)]$ and $[F_1 , {\pi_1}(e)]$
 are both rank-one operators. 
Hence 
$[F_1 , {\pi_1}(a)]$ 
is compact for all 
$a$ 
in the dense subalgebra of elements of the form
${\Sigma}_{n \in {\bf Z}} ({a_n}{S^n} +{b_n}{S^n}e)$,
where 
$\{ {a_n} \}, \{ {b_n} \} \in S({\bf Z})$, the
Schwarz functions on 
${\bf Z}$, and so for all $a \in {C^{*}}(\Gamma)$.  
 Hence  
${\bf w}_1 =(H, {\pi_1},F_1 , \gamma)$ 
is an even Fredholm module over 
$A$.

We calculate the pairing of the Chern character of ${\bf w}_1$ in even periodic cyclic cohomology with
the elements of 
${K_0}(A)$ 
represented by the projections 
$1$, 
${P_1} = {\frac{1}{2}}(1 + e)$,
${P_2} = {\frac{1}{2}}(1 + eS)$. 
We have 
\begin{equation}
{\pi_1}({P_1}) =
{\frac{1}{2}}
\left(
\begin{array}{cc}
1+e & 0 \cr
0 & 1-e \cr
\end{array}
\right), \quad
{\pi_1}({P_2}) =
{\frac{1}{2}}
\left(
\begin{array}{cc}
1+eS & 0 \cr
0 & 1-eS \cr
\end{array}
\right),
\end{equation}

\begin{equation}
[F_1 , {\pi_1}({P_1})]=
-{\frac{1}{2}}({F_0}e+e{F_0})
\left(
\begin{array}{cc}
0 & 1 \cr
1 & 0 \cr
\end{array}
\right)
= -i{P_{0,0}}
\left(
\begin{array}{cc}
0 & 1 \cr
1 & 0 \cr
\end{array}
\right),
\end{equation}
 where we denote by $P_{i,j}$ the rank one operator 
 $P_{i,j} e_m = \delta_{j,m} e_i$. 
Hence  
${[F_1 ,{\pi_1}({P_1})]}^{2k} $
$ = {(-1)}^{k} {P_{0,0}} {I_2} $,
so
\begin{equation}
Tr(
\gamma {\pi_1}({P_1}) {{[F , {\pi_1}({P_1})]}^{2k}})=
Tr(
{(-1)}^{k} {P_{0,0}}
\left(
\begin{array}{cc}
1 & 0 \cr
0 & 0 \cr
\end{array}
\right)
) = (-1)^k.
\end{equation}

Therefore
 $< {ch_{*}}( {\bf w}_1) , [{P_1}]> = 1$.
It follows that 
${\bf w}_{1}$ 
is a nontrivial element of 
${KK^0}(A, {\bf C})$.
The same calculation for 
${P_2}$ 
gives us that
\begin{equation}
[F_1 , {\pi_1}({P_2})] = 
-{\frac{1}{2}} i ({F_0}eS+ eS{F_0})
\left(
\begin{array}{cc}
0 & 1 \cr
1 & 0 \cr
\end{array}
\right)=0,
\end{equation}
since 
$FeS+eSF=0$.  
Hence 
$<{ch_{*}}({\bf w}_{1}),[ {P_2}]> = 0$.
Since
$[F_1 , {\pi_1}(1)] = 0$
it is obvious that 
$<{ch_{*}}({\bf w}_1),[1]> = 0$.

We obtain the third Fredholm module 
${\bf w}_{2}$ 
from ${\bf w}_1$  as follows.
 We pull back ${\bf w}_1$ via the automorphism 
$\alpha_{-1} \in Aut( \Gamma)$, defined by
 \begin{equation}
 \alpha_{-1} (S) = S,\quad
\alpha_{-1} (e) = {S^{-1}}e.
\end{equation}
 Then ${\bf w}_2 = {\alpha_{-1}^{*}} ( {\bf w}_1 )$ is as follows.
We again have  
$A$ 
acting on 
${l^2}({\bf Z})$, 
this time via
\begin{equation} 
S{e_n} = {e_{n+1}},\quad
e{e_n}={e_{-(n+1)}}
\end{equation} 
Take 
$H$, 
$F$ 
and 
$\gamma$ 
as before, 
and define the *-homomorphism 
${\pi_2} : A \rightarrow {\bf B}(H) $ 
by
\begin{equation}
{\pi_2}(S)=
\left(
\begin{array}{cc}
S & 0 \cr
0 & S \cr
\end{array}
\right), \quad
{\pi_2}(e) =
\left(
\begin{array}{cc}
e & 0 \cr
0 & -e \cr
\end{array}
\right),
\end{equation}
i.e. the same as 
$\pi_1$,
except the representation on 
${l^2}({\bf Z})$
is different. 
Then 
${\bf w}_{2}=$ 
 $(H, {\pi_2},F_2 = F_1, \gamma)$
 $ = (H, {\pi_1} \circ {\alpha_{-1}}, F_2, \gamma)$ 
 is an even Fredholm module over 
$A$,
 as can be verified exactly as before. 
To distinguish 
${\bf w}_{2}$
 from our previous example, 
we calculate the pairing of its Chern character with the projections
$1$, 
${P_1}$ and 
${P_2}$. 
 We note that 
\begin{equation}
\alpha_{-1} (1) =1,\quad
\alpha_{-1} (P_1) = \alpha_{-1}( {\frac{1}{2}}(1 + e)) = {\frac{1}{2}}(1 + {S^{-1}}e) \sim P_2,
\end{equation} 
\begin{equation}
\alpha_{-1}(P_2) =  \alpha_{-1}( {\frac{1}{2}}(1 + Se)) = {\frac{1}{2}}(1 + e) = P_1.
\end{equation}
 Hence 
\begin{equation}
<{ch_{*}}({\bf w}_2),[1]> = <{ch_{*}}({\bf w}_1),[\alpha_{-1}(1)] > =
<{ch_{*}}({\bf w}_1),[1]> =1,
\end{equation}
\begin{equation}
<{ch_{*}}({\bf w}_2),[{P_1}]> = <{ch_{*}}({\bf w}_1),[\alpha_{-1}({P_1})]> = <{ch_{*}}({\bf w}_2),[{P_2}]> = 0,
\end{equation}
\begin{equation}
<{ch_{*}}({\bf w}_2),[{P_2}]> = <{ch_{*}}({\bf w}_1),[\alpha_{-1}({P_2})]> = <{ch_{*}}({\bf w}_1),[{P_1}]> =1.
\end{equation}

We summarize the results of all these calculations in the following table.

\begin{Prop}
\label{dihedralkthykhompairings}
The pairings of the generating Fredholm modules with K-theory are given by :

$\begin{diagram}
 &  &1 	& & {\frac{1}{2}}(1+e)	 & & {\frac{1}{2}}(1+Se) \\
{\bf w}_0 & &1 & & 1 & & 1\\
{\bf w}_1 &  &0 & & 1 & & 0 \\
{\bf w}_2 &  &0 & & 0 & & 1 \\
\end{diagram}$
\end{Prop}

Since ${\bf w}_0$ 
pairs nontrivially with 
$[1] \in {K_0}(A)$, 
we can see that 
${\bf w}_0$, 
${\bf w}_1$ 
and 
${\bf w}_2$ 
in fact generate 
${KK^0}(A,{\bf C})$ 
as an abelian group - 
since their pairings with each of the projections are either 0 or 1, 
by Connes' index theorem \cite{connes94}, p296,
they are generators, 
rather than just generating a copy of 
${\bf Z}^3$ 
inside 
${KK^0}(A,{\bf C})$. 
 This completes the proof of Theorem \ref{khominfdih}.
\end{proof}

Note that under the inclusion map $i : C({\bf T}) \hookrightarrow A$, 
we have ${i^{*}}({\bf w}_0) = {\bf z}_0 \in {KK^0}(C({\bf T}),{\bf C})$, 
while 
\begin{equation}
{i^{*}}({\bf w}_1) = {i^{*}}({\bf w}_2) = {\bf 0} \in {KK^0}(C({\bf T}),{\bf C}).
\end{equation}

It would be highly desirable to generalize the methods of this section to more general crossed products by ${\bf Z}_2$. 
For a general algebra $A$, 
with an action $\sigma$ of ${\bf Z}_2$, it would be very useful to construct a map on K-homology
\begin{equation}
{KK^{*}}(A,{\bf C}) \rightarrow {KK^{*}}(A \times_{\sigma} {\bf Z}_2, {\bf C}).
\end{equation}
This can be done if $A$ is the group C*-algebra of a finitely generated abelian group, however  we have been unable to find any more general construction. 
This contrasts sharply  with the case of crossed products by ${\bf Z}$.

\section{Cyclic cohomology}

In this section we calculate the cyclic cohomology of the group ring 
${\bf C}\Gamma$ of $\Gamma$, 
and of a ``smooth subalgebra'' 
$A^{\infty}$ of the group C*-algebra.
 We define
\begin{equation}
\label{smoothdihedral}
{A^{\infty}} = \{ \Sigma( {a_m}{S^m} + {b_m}{S^m}e) \} 
\end{equation} 
 where 
$\{ {a_m} \}$ and $\{ {b_m} \}$ are Schwarz functions on ${\bf Z}$. 
 We note that for the obvious length function $L$ on $\Gamma$ defined by 
\begin{equation} 
 L(S^m) = |m|,\quad
 L({S^m}e) = |m| +1
\end{equation}
  $A^{\infty}$ is the space of rapidly decreasing functions on $\Gamma$ with repect to $L$. 
Since $A^{\infty}$ is a *-subalgebra of ${C_r^{*}}(\Gamma)$, 
it follows that  $\Gamma$ has property (RD) of Jolissaint \cite{jo}.
 Hence 
$A^{\infty}$ 
coincides with Jolissaint's algebra 
$H_{L}^{\infty}(\Gamma)$ 
and is therefore  a dense *-subalgebra of 
${C_r^{*}} (\Gamma)$ 
closed under holomorphic functional calculus. 
 We have 
\begin{equation}
{\bf C}\Gamma \subset A^{\infty} \subset {C^{*}}( \Gamma )
\end{equation}
with the inclusion of $A^{\infty}$ in ${C^{*}}(\Gamma)$ inducing an isomorphism on K-theory. 
Since $A = {C^{*}}(\Gamma)$ is a nuclear C*-algebra, its cyclic cohomology is given by 
\begin{equation}
\label{ccnuclear}
{HC^{n}}(A) \cong
\left\{ 
\begin{array}{cc}
{HC^0}(A) & : n \,even\\
0 & : n \,odd\\
\end{array}
\right.
\end{equation}
 with ${HC^0}(A)$ being generated by the traces on $A$. 
This remark is true for arbitrary nuclear C*-algebras \cite{connes85}, p132. 
The dense *-subalgebras ${\bf C}\Gamma$ and $A^{\infty}$ are not norm-closed and have much larger and more interesting cyclic cohomology. 

We begin by calculating the cyclic cohomology of ${\bf C}\Gamma$. 
We do this in two different ways. 
First of all via homological algebra, 
using Burghelea's theorem for cyclic homology and cohomology of group rings. 
Second, by a direct calculation of the cyclic 0-cocycles and 1-cocycles, and then a ``bare hands'' proof that every 1-cocycle is a 1-coboundary. 
 This approach explicitly gives us the generating cyclic cocycles. 
 We then investigate which of the cyclic cocycles on the group ring ${\bf C}\Gamma$ extend to  $A^{\infty}$, 
giving the cyclic cohomology 
${HC^{n}}( {A^{\infty}})$, for $n \geq 1$.

 Cyclic homology of group rings was calculated by Burghelea \cite{bu85} and stated for cyclic cohomology by Connes \cite{connes94}, p213. 
The result we use is as follows. 
Let $G$ be a countable discrete group. 
Given 
$g \in G$, 
let 
${C_g} = \{ h \in G : gh=hg \}$ 
be the centralizer of $g$. 
Let 
${N_g} = {C_g}/ {g^{\bf Z}}$ 
be the quotient of 
${C_g}$ by the central normal subgroup generated by $g$. 
Now let $<G>$ be the set of conjugacy classes of elements of $G$, 
let $<G>'$ be the set of conjugacy classes of elements of finite order, 
and let $<G>''$ be the set of conjugacy classes of elements of infinite order. 

\begin{Thm}[Burghelea] \cite{connes94}
\label{burg}
 The cyclic cohomology of ${\bf C}G$ is given by 
\begin{equation}
{HC^{*}}({\bf C}G) = {\Pi_{ {\hat{g}} \in <G>'} ( {H^{*}}( {N_g};{\bf C}) \otimes {HC^{*}}({\bf C}) )} \times  {\Pi_{ {\hat{g}} \in <G>''} {H^{*}}( {N_g}; {\bf C})}
\end{equation}
\end{Thm} 

Here ${H^{*}}(G;{\bf C})$ denotes the cohomology of the group $G$ with coefficients in ${\bf C}$ (as a trivial $G$-module).

We apply this to the situation $G = \Gamma$. 
As before we write the generators of $\Gamma$  as $S$ and $e$ with the relations
\begin{equation}
e^2=1,\quad 
 eSe = S^{-1}.
\end{equation}
 Let $<\Gamma>$ be  the set of all conjugacy classes.
 Since 
$e \sim {S^{2n}}e$, 
$Se \sim {S^{2n+1}}e$ 
for all $n$, 
and 
${S^{n}} \sim {S^{-n}}$, 
 the conjugacy classes  
$[1]$, $[e]$ and  $[Se]$ 
correspond to elements of finite order, 
and 
${ \{ [{S^n}] \} }_{n \geq 1}$ 
correspond to elements of infinite order. 
  We have  $<\Gamma>' = \{ [1],[e],[Se] \}$. 
Obviously  $N_1 = \Gamma$,  
 while $N_{e} = {N_{Se}} = {1}$ (the trivial group), and     
\begin{equation}
{H^n}( 1 ; {\bf C} ) = 
\left\{ 
\begin{array}{cc}
{\bf C} & : {n=0} \\
0 &  :  {n \geq 1} \\
\end{array}
\right. 
\end{equation}
 We need the following :

\begin{Prop} 
${H^n}( \Gamma ; {\bf C} ) = 
\left\{ 
\begin{array}{cc}
{\bf C} & : {n=0} \\
0 &  :  {n \geq 1} \\
\end{array}
\right.$
\end{Prop} 
 \begin{proof} We know \cite{weibel}, p170 that for finite groups $G$ and  $H$, and any left $G*H$-module $A$, that for $n \geq 2$ the group cohomology of the free product $G*H$ is given by  
\begin{equation}
{H^n}( G*H ; A) = {H^n}(G;A) \oplus {H^n}(H;A)
\end{equation} 
 Since $\Gamma = {\bf Z}_2 * {\bf Z}_2 $, and 
\begin{equation}
{H^n}( {\bf Z}_2  ; {\bf C} ) = 
\left\{ 
\begin{array}{cc}
{\bf C} & : {n=0} \\
0 & : {n \geq 1} \\
\end{array}
\right.
\end{equation}
 it follows that ${H^n}(\Gamma ; {\bf C}) = 0$ for all $n \geq 2$. 
By definition, ${H^0}(\Gamma ; {\bf C} )= {\bf C}^{\Gamma} = {\bf C}$ since ${\bf C}$ is a trivial $\Gamma$-module. 
 To find ${H^1}(\Gamma ; {\bf C})$ we use the Hochschild-Serre spectral sequence. 
Recall that, for an exact sequence of groups 
\begin{equation}
0 \rightarrow H \rightarrow G \rightarrow G/H \rightarrow 0 
\end{equation}
 and some $G$-module $M$, 
${H^q}( H;M)$ is a $G/H$-module, and 
${E_2^{p,q}} = {H^p}( G/H ; {H^q}(H;M) ) $ converges to ${H^{p+q}}(G;M)$. 
 In our case we have an exact sequence 
\begin{equation}
0 \rightarrow {\bf Z} \rightarrow \Gamma \rightarrow {\bf Z}_2 \rightarrow 0 
\end{equation}
 so 
\begin{equation} 
{E_2^{p,q}} = {H^p}( {\bf Z}_2 ; {H^q}( {\bf Z};{\bf C}))
\end{equation}
 and 
\begin{equation} 
{H^q}( {\bf Z}; {\bf C} ) = 
\left\{ 
\begin{array}{cc}
{\bf C} & : {q=0,1} \\
0 &  : {q \geq 2} \\
\end{array}
\right.
\end{equation}
 So ${E_2^{p,q}}=0$ for $q > 1$. 
 Now, ${E_2^{1,1}} ={H^1}( {\bf Z}_2 ; {H^1}({\bf Z};{\bf C}))$ 
with the action of ${\bf Z}_2$ on ${H^1}({\bf Z};{\bf C}) \cong {\bf C}$ being by inner automorphisms (inversion) so 
${E_2^{1,1}} = 0$. 
 Finally, 
${E_2^{0,0}} = {H^0}({\bf Z}_2 ; {\bf C}) = {\bf C}$ is the only non-zero ${E_2^{p,q}}$. 
So everything is over at the first step of the spectral sequence.  
\end{proof}

 We know \cite{connes94}, p192, that the cyclic cohomology of the complex numbers ${\bf C}$ is given by: 
\begin{equation}
{HC^n}({\bf C} ) = 
\left\{ 
\begin{array}{cc}
{\bf C} & : {n=0,2,4..} \\
0 &  :  {n =1,3,5..} \\
\end{array}
\right.      
 \end{equation}
Hence 
\begin{equation}
{\Pi_{ {\hat{g}} \in <G>'} ( {H^{*}}( {N_g};{\bf C}) \otimes {HC^{*}}({\bf C}) )} =
\left\{ 
\begin{array}{cc}
{\bf C}^3   & : {n=0,2,4..} \\
0 &  :  {n =1,3,5..} \\
\end{array} 
\right.
 \end{equation} 
We also have $<G>'' = {\{ [S^m] \}}_{m \geq 1}$. 
Now, $N_{S^m} \cong {\bf Z}_m$ ($m \geq 1$), and 
\begin{equation}
{H^n}( {\bf Z}_m  ; {\bf C} ) = 
\left\{ 
\begin{array}{cc}
{\bf C} & : {n=0} \\
0 &  :  {n \geq 1} \\
\end{array}
\right.
\end{equation}
 so 
\begin{equation}
{\Pi_{ {\hat{g}} \in <G>''} {H^{*}}( {N_g}; {\bf C})} =
\left\{
\begin{array}{cc}
{\Pi_{m \geq 1}} {\bf C} & : {n=0} \\
0 &  :  {n \geq 1} \\
\end{array}
\right.
\end{equation}
 We have therefore proved the following : 
\begin{Prop}
\label{ccCGamma} 
The cyclic cohomology of the group ring ${\bf C}\Gamma$ is given by 
\begin{equation}
{HC^n}({\bf C}\Gamma ) = 
\left\{
\begin{array}{cc}
{\bf C}^3 \times {\Pi_{p \geq 1}} {\bf C}  & : {n=0} \\
0 &  : {n =1,3,5..} \\
{\bf C}^3 & : {n=2,4,6..}\\
\end{array}
\right.
\end{equation}
\end{Prop}

The reason for the above notation for ${HC^0}({\bf C}\Gamma)$ is that we distinguish three 0-cocycles ${\psi_0}$, $\psi_1$ and $\psi_2$  that we exhibit below. 
These cocycles generate the even periodic cyclic cohomology (as a ${\bf C}$-module.)
 The other generating 0-cocycles do not contribute to periodic cyclic cohomology (Prop \ref{ccgpring1}, Lemma \ref{sphikiszero}).

\begin{Cor}
The periodic cyclic cohomology of ${\bf C}\Gamma$ is given by 
\begin{equation}
{HC^{even}}( {\bf C}\Gamma) \cong {\bf C}^3 ,\quad
{HC^{odd}}( {\bf C}\Gamma) \cong {\bf 0}. 
\end{equation}
\end{Cor}

We now give a ``bare hands'' proof of Prop \ref{ccCGamma}.
 This direct calculation explicitly gives the 0-cocycles generating the even periodic cyclic cohomology. 

Let $\psi$ be a 
0-cocycle on ${\bf C}\Gamma$. Then $\psi(xy) = \psi(yx)$ for all $x$ and  $y$, hence $\psi$ is constant on conjugacy classes of $\Gamma$. 
 Taking $x=e$ and $y={S^n}e$, we see that $\psi(S^n) = \psi(S^{-n})$ for all $n$. 
Define $\psi(S^n) = a_n$, for some numbers $\{ a_n \}$ with $a_{-n} =  a_n$. 
Taking $x={S^m}$, $y={S^n}e$ we see that 
$\psi({S^{m+n}}e) = \psi({S^{n-m}}e)$ for all $m$, $n$. 
Hence  
$\psi({S^{m+2n}}e) = \psi({S^m}e) $ for all $n$. 
So take 
$\psi(e) = b_0$, $\psi(Se) = b_1$. 
 The pairings of $\psi$ with the projections $P_0 =1$, $P_1 = {\frac{1}{2}}(1+e)$, ${P_2} = {\frac{1}{2}}(1 + Se)$ are given by 
$\psi(1) = a_0$, $\psi( {\frac{1}{2}}(1+e)) = {\frac{1}{2}}(a_0 + b_0 )$, and 
 $\psi( {\frac{1}{2}}(1+Se)) = {\frac{1}{2}}(a_0 + b_1 )$.
 We denote by $\psi_0$, $\psi_1$ and  $\psi_2$ the cyclic 0-cocycles corresponding to $a_0 =1$, $b_0 =2$ and $b_1 = 2$ respectively (in each case all the other $a_n$ and $b_n$ are taken to be zero). 
 Then ${\psi_i}( P_j ) = \delta_{i,j}$. 
We will see later that these distinguished cyclic 0-cocycles generate the even periodic cyclic cohomology.  
We have therefore shown directly that 
\begin{equation}
{HC^0}({\bf C}\Gamma) \cong {\bf C}^3 \times {\Pi_{p \geq 1}} {\bf C}.
 \end{equation}

Now we calculate the 1-coboundaries. 
Let  $\psi$ be a linear functional (not assumed to be a cocycle) 
$\psi : {\bf C} \Gamma \rightarrow {\bf C}$. 
Then $\psi$ is completely determined by two (not necessarily bounded) sequences $\{ a_m \}$, $\{ b_m \}$ ($ m \in {\bf Z}$) with  
\begin{equation} 
\psi(S^m) = a_m , \quad \psi({S^m}e) = b_m.
\end{equation}  
Now $b \psi$ is the 1-coboundary given by $b \psi (x,y) = \psi(xy) - \psi(yx)$, hence we see that 
 \begin{equation}
b \psi(S^m, S^n) = \psi(S^{m+n}) - \psi(S^{m+n}) = 0,
\end{equation}
\begin{equation}
b \psi(S^m, {S^n}e) = \psi({S^{m+n}}e) - \psi({S^{n-m}}e) = {b_{m+n}}-{b_{n-m}},
\end{equation}
\begin{equation}
b \psi({S^m}e, {S^n}e) = \psi(S^{m-n}) - \psi(S^{n-m}) = {a_{m-n}}-{a_{n-m}}.
\end{equation}
 Now let $\phi$ be a 1-cocycle. By definition  
\begin{equation}
\phi(x,y) = - \phi(y,x), \quad 
\phi(xy,z) - \phi(x,yz) + \phi(zx,y) = 0. 
\end{equation}

{\bf Case 0:} $x={S^p}$, $y={S^q}$, $z={S^r}$. 
 We need $\phi$ to restrict to a cyclic 1-cocycle on the copy of ${\bf C}{\bf Z}$ generated by $S$. Hence
$\phi ( {S^m}, {S^n} ) = k n {\delta_{m+n,0}}$ 
 for some $k \in {\bf C}$. 

{\bf Case 1:} 
 For $x={S^p}$, $y={S^q}$, $z={S^r}e$, we have   
 \begin{equation}
\phi( {S^{p+q}}, {S^r}e ) - \phi( {S^p}, {S^{q+r}}e) + \phi( {S^{r-p}}e, {S^q})
 =0.
 \end{equation}
 Writing $f(m,n) = \phi( {S^m}, {S^n}e)$ , and using the cyclic property, we find that 
 \begin{equation}
 f(p+q,r) -f(p, q+r) -f(q,r-p) =0.
\end{equation}
 Setting $q=0$ gives 
$f(p,r) - f(p,r) = 0 = f(0, r-p)$ 
for all $r$, $p$, i.e. $f(0,m)=0$ for all $m$. 
 Setting $q=1$ gives  
$f(p+1,r) - f(p,r+1) = f(1,r-p)$, 
 and we can solve this relation to get 
 \begin{equation}
 f(p+1,r) = {\Sigma_{j=0}^{p}} f(1,r-p+2j).
 \end{equation} 
 Writing $\phi(S, {S^m}e) = c_m = f(1,m)$, we have 
 \begin{equation}
 \phi({S^m}, {S^n}e) = c_{n+m-1} + c_{n+m-3} + c_{n+m-5} + .. + c_{n-m+1}.
 \end{equation}

{\bf Case 2:} $x={S^p}$, $y={S^q}e$, $z={S^r}e$. 
This gives 
\begin{equation}
\phi({S^{p+q}}e, {S^r}e) - \phi({S^p},{S^{q-r}}) + \phi({S^{r-p}}e, {S^q}e) =0.
 \end{equation}
 For $p+q-r =0$, this gives 
\begin{equation}
\phi({S^r}e, {S^r}e) - \phi({S^p},{S^p}) + \phi({S^q}e,{S^q}e) = -\phi({S^p},{S^p}) = -kp = 0,
\end{equation}
 hence we must have $k=0$, i.e. $\phi( {S^m},{S^n}) = 0$ for all $m$, $n$. 

 For $p+q-r \neq 0$, we have 
$\phi({S^{p+q}}e, {S^r}e) = \phi({S^q}e, {S^{r-p}}e)$. 
Writing $g(m,n) = \phi({S^m}e, {S^n}e)$ gives  
\begin{equation} 
g(m,n) = - g(n,m), \, g(q+p,r) = g(q,r-p),
\end{equation}
so 
\begin{equation} 
g(q+1,r) = g(q,r-1) = ... = g(0,r-q-1).
\end{equation}
 Taking $g(0,n) = {d_n}$, with $d_{-n} = - d_n$, we have  
\begin{equation} 
\phi ({S^m}e, {S^n}e ) = d_{n-m}. 
\end{equation}
 Since $d_0 =0$ this is true for all $m$, $n$. 

{\bf Case 3:} $x={S^p}e$, $y={S^q}e$, $z={S^r}e$. This reduces to Case 1. 

Hence any cyclic 1-cocycle $\phi$ is uniquely determined by two sequences $\{ {c_n} \}$, $\{ {d_n} \}$ with $d_{-n} = - d_n$ and is given by 
\begin{equation}
  \phi( {S^m},{S^n}) = 0, \quad
 \phi({S^m}, {S^n}e) = {\Sigma_{k=0}^{m-1}} c_{n+m-1-2k}, \quad 
 \phi ({S^m}e, {S^n}e ) = d_{n-m}.
 \end{equation}

\begin{Prop} 
Every 1-cocycle is a coboundary, so ${HC^1}({\bf C}\Gamma) =0$.
\end{Prop} 
\begin{proof}
Given a 1-cocycle $\phi$, we show that we always have $\phi = b \psi$, for some linear map $\psi : {\bf C}\Gamma \rightarrow {\bf C}$, with 
$\psi(S^m) = a_m$,  $\psi({S^m}e) = b_m$. 
 By our previous work we have 
\begin{equation} 
\phi( {S^m}, {S^n} ) = 0 = b \psi ( {S^m}, {S^n}), 
 \end{equation}
\begin{equation} 
\phi ( {S^m}, {S^n}e ) = 
 {\Sigma_{k=0}^{m-1}} {c_{n-m+1+2k}}, 
 \end{equation}
\begin{equation} 
\phi ( {S^m}e, {S^n}e ) = {d_{n-m}}. 
 \end{equation}
Now define 
\begin{equation}
a_m = 
\left\{ 
\begin{array}{cc}
-{d_m} &  : m>0 \\
0 &  :  m \leq 0 \\
\end{array}
\right.   
\end{equation} 
\begin{equation} 
 b_{2m} = 
\left\{ 
\begin{array}{cc}
c_1 + c_3 + ... + c_{2m-1}  & :  m>0 \\
0 & : m=0\\
-c_{-1} - c_{-3} - ... - c_{2m+1} & :  m < 0 \\
\end{array}
\right.   
\end{equation}

\begin{equation} 
 b_{2m+1} = 
\left\{ 
\begin{array}{cc}
c_0 + c_2  + ... + c_{2m}  & :  m>0 \\
-c_{-1} - c_{-3} - ... - c_{2m+2} & :  m \leq 0 \\
\end{array}
\right.   
\end{equation}
 So we have shown that every 1-cocycle is in fact a coboundary, so ${HC^1}({\bf C}\Gamma)=0$. 
\end{proof}

\begin{Prop} 
\label{ccgpring1}
 The cyclic 0-cocycles $\psi_0$, $\psi_1$, $\psi_2$ generate the even periodic cyclic cohomology ${HC^{even}}( {\bf C}\Gamma ) \cong {\bf C}^3$.
\end{Prop}
 \begin{proof}  From our previous work we see immediately that 
 ${HC^{even}}( {\bf C}\Gamma ) \cong {\bf C}^3$, 
while\\
  ${HC^{odd}}( {\bf C}\Gamma ) \cong 0$. 
 The three distinguished 0-cocycles are ``dual'' to the generators of\\ 
 ${K_0}( {C_r^{*}}(\Gamma))$, 
 and Connes' periodicity operator 
 $S : {HC^i}(A) \rightarrow {HC^{i+2}}(A)$ 
maps them to cohomologically nontrivial even cocycles of higher degree.  
We see that for each $n >0$ we have  
${HC^{2n}}({\bf C}\Gamma) \cong {\bf C}^3$, 
with generators ${S^n}{\psi_0}$, ${S^n}{\psi_1}$ and  ${S^n}{\psi_2}$. 
\end{proof}

 Recall that we found that the zeroth cyclic cohomology was infinite dimensional. 
We now show directly that the generating 0-cocycles not in the span of the three given above are mapped to 
$0 \in {HC^2}({\bf C}\Gamma)$ by $S$. 
 Denote by ${\psi_k}$ the cyclic 0-cocycle corresponding to $a_k = 1$, 
all others zero (we take $k \neq 0$). 
Explicitly, 
\begin{equation}
{\psi_k}(S^m) = {\delta_{k,m}}+{\delta_{k,-m}}, \,{\psi_k}(e) = 0 = {\psi_k}(Se)
\end{equation}
\begin{Lemma}
\label{sphikiszero}
 The 2-cocycle  $S{\psi_k} \in {HC^2}({\bf C}\Gamma)$ is a coboundary, for $k \neq 0$. 
\end{Lemma}
\begin{proof} 
 By definition $S{\psi_k} (x,y,z) = {\psi_k}(xyz)$. 
Let $\phi : {\bf C}\Gamma \times {\bf C}\Gamma \rightarrow {\bf C}$ be a linear functional, with $\phi(x,y) = - \phi(y,x)$.  
Suppose that 
\begin{equation} 
\phi (S^m, S^n) = \alpha_{m,n}, \quad
\phi (S^m, {S^n}e) = \beta_{m,n}, \quad
\phi ({S^m}e, {S^n}e) = \gamma_{m,n}.
\end{equation}
 Obviously 
$\alpha_{n,m} = - \alpha_{m,n}$, 
and $\gamma_{n,m} = - \gamma_{m,n}$. 
We want to solve 
\begin{equation} 
 b\phi (x,y,z) = \phi(xy,z) - \phi(x,yz) + \phi(zx,y) = S{\psi_k}(x,y,z) = {\psi_k}(xyz).
 \end{equation}

{\bf Case 0:} 
$x= {S^p}$, $y={S^q}$, $z={S^r}$. This gives 
\begin{equation} 
 \alpha_{p+q,r} - \alpha_{p,q+r} + \alpha_{p+r,q} = {\psi_k}({S^{p+q+r}}) = 
\left\{ 
\begin{array}{cc}
1 &    : p+q+r = \pm k \\
0 &  : otherwise \\
\end{array}
\right.    
\end{equation} 
which has the solution 
\begin{equation}
\alpha_{m,n} = 
\left\{ 
\begin{array}{cc}
{\frac{m-n}{m+n}}  &  :  m+n = \pm k  \\
0 &  :  otherwise \\
\end{array}
\right.   
\end{equation}

{\bf Case 1:}
$x= {S^p}$, $y={S^q}$, $z={S^r}e$. This gives 
\begin{equation}
\beta_{p+q,r} - \beta_{p,q+r} - \beta_{q,r-p} = 
{\psi_n}( {S^{p+q+r}}e) = 0,
\end{equation}
 which has the solution $\beta_{m,n} = 0$, for all $m$, $n$. 

{\bf Case 2:}
$x= {S^p}$, $y={S^q}e$, $z={S^r}e$. This gives 
\begin{equation} 
\gamma_{p+q,r} - \alpha_{p,q-r} - \gamma_{q,r-p} = a_{p+q-r}.
\end{equation}
Hence
\begin{equation} 
\gamma_{p+q,r}  - \gamma_{q,r-p} = 
\left\{ 
\begin{array}{cc}
{\frac{2p}{p+q-r}}  &   :  p+q-r = \pm k  \\
0 & : otherwise \\
\end{array}
\right. 
\end{equation}
 This has the solution 
\begin{equation}
\gamma_{m,n} = 
\left\{ 
\begin{array}{cc}
{\frac{2n}{k}} - {c_k}   & :  m-n = k  \\
{\frac{-2m}{k}} + {c_k}   &  :  m-n = -k \\
0 & : otherwise \\
\end{array}
\right. 
\end{equation}
 where $c_k$ is a scalar. 
Hence we have shown that for any of the 0-cocycles ${\psi_k}$, with $k \neq 0$, the corresponding 2-cocycle $S{\psi_k}$ is a coboundary. 
So for each such $k$, and any $n \geq 1$, 
${S^n}{\psi_k} = 0 \in {HC^{2n}}( {\bf C} \Gamma )$. 
\end{proof}

Now we consider the cyclic cohomology of the ``smooth subalgebra''  $A^{\infty}$ of the group C*-algebra we defined earlier (\ref{smoothdihedral}). 
Obviously any cyclic cocycle on $A^{\infty}$ restricts to a cyclic cocycle on ${\bf C} \Gamma$, 
 and the cyclic 0-cocycles $\psi_0$, $\psi_1$ and $\psi_2$ we exhibited generating the periodic even cyclic cohomology of  
${\bf C}\Gamma$ all extend to 0-cocycles on $A^{\infty}$, via
 \begin{equation}
\psi( \Sigma  \alpha_p S^p + \beta_p S^p e) = 
\Sigma  \alpha_p \psi(S^p) + \beta_p \psi(S^p e) .
\end{equation} 
 However, this is not necessarily enough to determine the cyclic cohomology of $A^{\infty}$. 
 A cyclic cocycle on $A^{\infty}$ is not necessarily uniquely determined by its restriction to ${\bf C} \Gamma$. 
 There is the possibility of nonzero ``ghost'' cocycles that vanish identically on the group ring. 
 Since the three projections that generate $K_0 (A)$ are all elements of ${\bf C} \Gamma$, such ghosts will not be detectible by means of pairing with K-theory. 

Nevertheless, we conjecture that:

\begin{Conj} 
\label{ccsmoothdih}
The cyclic cohomology of $A^{\infty}$ is given by 
\begin{equation}
{HC^{n}}( A^{\infty} ) = 
\left\{ 
\begin{array}{cc}
{\bf C}^3 & : n \,\,even, n \geq 2\\
0 & : n \,\,odd, n \geq 1 \\
\end{array}
\right. 
\end{equation} 
\end{Conj}

It is clear that ${HC^0}(A^{\infty})$ is very large. 
Any 0-cocycle $\psi$ on ${\bf C}\Gamma$ is uniquely determined by a (not necessarily bounded) sequence $\{ {a_n} \}$, 
with $a_{-n} = - a_n$, and scalars $b_0$, $b_1$. 
If the $a_n$ are of polynomial growth in $n$, it is clear that $\psi$ extends to a well-defined cyclic cocycle  on $A^{\infty}$, otherwise not.

\section{${C^{*}}( {\bf Z} \times_{\sigma} {\bf Z})$, a noncommutative orbifold}

Following on from our work on the infinite dihedral group, 
 in this section we study the group $G = {\bf Z} \times_{\sigma} {\bf Z}$.

Since $G$ is amenable, the full and reduced group C*-algebras are isomorphic. 
The group C*-algebra $B = {C^{*}}(G)$ is generated by unitaries $U$ and $V$ satisfying $VU = {U^{*}}V$. 
Let  $A = $ ${C^{*}}(\Gamma)$ be the group C*-algebra of the infinite dihedral group, generated by unitaries $S$ and $e$ with $e={e^{*}}$, $e^2 =1$ and $eSe = S^{*}$.
 We have a quotient *-homomorphism 
$ B \rightarrow A$ given by 
$ U \mapsto S$, $V \mapsto e$. 

\section{K-theory}

The group C*-algebra $B$ is a crossed product $C({\bf T}) \times_{\sigma} {\bf Z}$. Therefore we can use the Pimsner-Voiculescu sequence to calculate the K-theory. 

\begin{diagram}
{K_0} ( C( {\bf T} )) & \rTo^{ id - {\sigma_{*}}} & {K_0} (C( {\bf T} )) & \rTo^{ i_{*}} 
& {K_0} (B) \\
\uTo^{\delta_1} &    &  &  & \dTo^{\delta_0} \\
{K_1} (B) & \lTo^{ i_{*}} & {K_1} (C( {\bf T})) & \lTo^{id - {\sigma_{*}}} & {K_1} (C({\bf T})) 
\end{diagram}

We consider $C({\bf T})$ to be the C*-subalgebra of $B$ generated by the unitary $U$, and $i$ is the inclusion map 
$i : C({\bf T}) \hookrightarrow B$ given by $U \mapsto U$. 
The action $\sigma$ of ${\bf Z}$ on $C({\bf T})$ is defined by 
$\sigma(f) (z) = f( \overline{z})$, or alternatively 
$\sigma(U) = {U^{*}}$. 
We know that ${K_0}( C({\bf T}) ) \cong {\bf Z}$, generated by $[1]$, and 
${K_1}( C({\bf T}) ) \cong {\bf Z}$, generated by $[U]$. 

We have ${\sigma_{*}}[1] = [ \sigma(1) ] = [1]$, hence 
the map $(id - {\sigma_{*}}) : {K_0}( C({\bf T}) ) \rightarrow {K_0}( C({\bf T}) )$ is the zero map. 
 Further, we have ${\sigma_{*}}[U] = [ \sigma(U) ] = [U^{*}] = -[U]$,
so the map $(id - {\sigma_{*}}) : {K_1}( C({\bf T}) ) \rightarrow {K_1}( C({\bf T}) )$ is the map 
${\bf Z} \rightarrow {\bf Z}$ given by $n \mapsto 2n$. 
So $\delta_0 : {K_0}(B) \rightarrow {K_1}(C({\bf T}))$ is the zero map, and hence 
${K_0}(B) \cong {\bf Z}$, generated by $[1]$. 
We will later confirm that $[1]$ is a nonzero generator of ${K_0}(B)$ by showing that $[1]$ pairs to 1 with an even Fredholm module over $B$. 

Further, we see that ${i_{*}}[U] = [U]$ is a nonzero torsion element of ${K_1}(B)$, with $[U]+[U]=0 \in {K_1}(B)$. 
Finally, the map $\delta_1 : {K_1}(B) \rightarrow {K_0}( C({\bf T}))$ is surjective.  
 So $[V]$ is a  nontrivial nontorsion generator of ${K_1}(B)$, 
with $\delta_1 [V] = \pm [1] \in {K_0}(C({\bf T}))$.  
We summarize these results in: 

\begin{Prop}
\label{kthyB}
${K_0}(B) \cong {\bf Z}$, generated by $[1]$, and 
${K_1}(B) \cong {\bf Z} \oplus {\bf Z}_2$, generated by $[V]$ and $[U]$, with $[U] + [U] =0$. 
\end{Prop}

\section{K-homology}

Since the K-groups have torsion, the universal coefficient theorem (Prop \ref{kgpsfreeab}) does not give us the K-homology for free. 
Instead, we consider the six term cyclic sequence on K-homology dual to the Pimsner-Voiculescu sequence on K-theory. 

\begin{diagram}
{KK^0}( C({\bf T}) , {\bf C}) & \lTo^{id - {\sigma}^{*}} & {{KK^0}( C({\bf T}), {\bf C})} & \lTo^{i^{*}} & {{KK^0}(B , {\bf C})}  \\
\dTo^{\partial_0} & & & & \uTo^{\partial_1} \\
{KK^1}( B, {\bf C}) & \rTo^{i^{*}} & {{KK^1}( C({\bf T}), {\bf C})} & \rTo^{id - {\sigma}^{*}} & {{KK^1}( C({\bf T}), {\bf C} )}\\
\end{diagram}

In Lemma \ref{onetorus}, we saw that both the even and odd K-homology of $C({\bf T})$ are isomorphic to ${\bf Z}$, with generators 
${\bf z}_0$ and ${\bf z}_1$ respectively.

We note straightaway the canonical even Fredholm module (Example \ref{evenfred}) 
${\bf w}_0 \in {KK^0}(B , {\bf C})$ corresponding to the *-homomorphism 
$\phi : B \rightarrow {\bf C}$ defined by $U, V \mapsto 1$.
We have
\begin{equation}
{\bf w}_0 = (H={\bf C}^2, \pi = \phi \oplus 0, 
F= \left(
\begin{array}{cc}
0 & 1 \cr
1 & 0 \cr
\end{array}
\right),
\gamma = 
\left(
\begin{array}{cc}
1 & 0 \cr
0 & -1 \cr
\end{array}
\right)
)
\end{equation}

\begin{Lemma} We have 
${i^{*}}({\bf w}_0) = {\bf z}_0$, and  
${\sigma^{*}}( {\bf z}_0 ) = {\bf z}_0$.
\end{Lemma}
\begin{proof}
Both these statements are immediate, the second because $\phi \circ \sigma = \phi$. 
\end{proof}

 It follows that 
$(id - {\sigma}^{*}) : {KK^0}( C({\bf T}) , {\bf C}) \rightarrow {KK^0}( C({\bf T}) , {\bf C})$ 
is the zero map. 
As described in Prop. \ref{partialzero}, 
 under the map $\partial_0$, we obtain a Fredholm module 
${\bf w}_1 = {\partial_0}( {\bf z}_0 ) \in {{KK^1}( B, {\bf C})}$, 
given by 
${\bf w}_1 = ({l^2}({\bf Z}), \pi_1 ', F)$, 
where 
${\pi_1 '}(V) e_n = e_{n+1}$, 
${\pi_1 '}(U) = I$, 
and $F{e_n} = sign(n) e_n$,  as before. 

\begin{Lemma}
 $i^{*} ( {\bf w}_1 ) = 0 \in {KK^1}( C({\bf T}), {\bf C})$.
\end{Lemma}
\begin{proof}
 We see that  
$i^{*} ( {\bf w}_1 ) = ({l^2}({\bf Z}), {\pi_1 '} \circ i, F)$ 
is a degenerate Fredholm module, since ${\pi_1 '} \circ i (U) = I$,  
and hence is zero in 
${KK^1}( C({\bf T}), {\bf C})$.
\end{proof}

\begin{Lemma}
We have $<{ch_{*}}( {\bf w}_0 ),[1]> =1$, and 
$<{ch_{*}}( {\bf w}_1 ),[V]> =1$, so 
${\bf w}_0$ and ${\bf w}_1$ are nonzero nontorsion generators of K-homology, and $[1]$ and $[V]$ are nonzero nontorsion generators of K-theory. 
\end{Lemma}
\begin{proof}
 This follows immediately from Examples \ref{evenfred} and \ref{oddfred}.
\end{proof}

This in fact completes the proof of Prop \ref{kthyB}. 
%\end{proof}

Finally,  
$\partial_1({\bf z}_1)$ is the Fredholm module 
\begin{equation}
\partial_1({\bf z}_1) = 
(
H = {l^2}({\bf Z}^2) \oplus {l^2}({\bf Z}^2), \pi \oplus \pi, 
\tilde{F_0} = 
\left(
\begin{array}{cc}
0 & {F_0} \cr
{F_0}^{*} & 0 \cr
\end{array}
\right),
\gamma =
\left(
\begin{array}{cc}
1 & 0 \cr
0 & -1 \cr
\end{array}
\right)
)
\end{equation}
 where 
$\pi(V) e_{p,q} = e_{p+1,q}$, 
$\pi(U) e_{p,q} = e_{p, q+ {(-1)}^p}$, and 
\begin{equation}
{F_0} e_{p,q} = 
\left\{
\begin{array}{cc}
{\frac{p+iq}{(p^2 + q^2)^{1/2}}} e_{p,q} & : (p,q) \neq (0,0) \cr
e_{0,0} & : (p,q)=(0,0) \cr
\end{array}
\right.
\end{equation} 

\begin{Lemma}
${i^{*}}( \partial_1({\bf z}_1) ) = {\bf 0} \in {KK^0}( C({\bf T}), {\bf C})$. 
\end{Lemma}
\begin{proof}
 We have 
\begin{equation}
{i^{*}}( \partial_1({\bf z}_1) ) 
= 
(
H = {l^2}({\bf Z}^2) \oplus {l^2}({\bf Z}^2), \pi \oplus \pi, 
\tilde{F_0} = 
\left(
\begin{array}{cc}
0 & {F_0} \cr
{F_0}^{*} & 0 \cr
\end{array}
\right),
\gamma =
\left(
\begin{array}{cc}
1 & 0 \cr
0 & -1 \cr
\end{array}
\right)
)
\end{equation}
 where $\pi(U) e_{p,q} = e_{p, q+ {(-1)}^p}$, and ${F_0}$ is as above. 
We construct a homotopy of Fredholm modules from 
${i^{*}}( \partial_1({\bf z}_1) )$ to a degenerate Fredholm module. 
For $0 \leq t \leq 1$, we define 
 \begin{equation}
{\bf y}_t
= 
(
H = {l^2}({\bf Z}^2) \oplus {l^2}({\bf Z}^2), \pi \oplus \pi, 
\tilde{F_t} = 
\left(
\begin{array}{cc}
0 & {F_t} \cr
{F_t}^{*} & 0 \cr
\end{array}
\right),
\gamma =
\left(
\begin{array}{cc}
1 & 0 \cr
0 & -1 \cr
\end{array}
\right)
)
\end{equation}
 where $F_t$ is given by 
\begin{equation}
{F_t} e_{p,q} = 
\left\{
\begin{array}{cc}
sign(p) e_{p,0} &: q=0 \cr
{\frac{p+i(1-t)q}{(p^2 + (1-t)^2 q^2)^{1/2}}} e_{p,q} & : q \neq 0 \cr
\end{array}
\right.
\end{equation} 
 Then ${\bf y}_0 =$ ${i^{*}}( \partial_1({\bf z}_1) )$, and 
${\bf y}_1$ is a degenerate Fredholm module, since $[F_1, \pi(U)]=0$. 
 Hence ${i^{*}}( \partial_1({\bf z}_1) )$ 
is a trivial element of K-homology.
\end{proof}

\begin{Prop}
\label{khomologyB}
We have ${KK^0}(B,{\bf C}) \cong {\bf Z} \oplus {\bf Z}_2$, generated by ${\bf w}_0$ and $\partial_1 ( {\bf z}_1 )$, with $\partial_1 ( {\bf z}_1 ) + \partial_1 ( {\bf z}_1 ) =0$, and 
${KK^1}(B,{\bf C}) \cong {\bf Z}$, generated by ${\bf w}_1$.
\end{Prop}

\section{Acknowledgements}
 I would like to thank my advisor Professor Marc Rieffel for his advice and support throughout my time at Berkeley. I am extremely grateful for his help. 
 I would also like to thank Erik Guentner, Nate Brown and Frederic Latremoliere for many useful discussions.

%\bibliographystyle{amsalpha}
%\bibliography{bibdata}
%\providecommand{\bysame}{\leavevmode\hbox to3em{\hrulefill}\thinspace}

\end{document}